\documentstyle{amsppt}
\topmatter
\title Isoperimetric Functions for Graph Products
\endtitle
\author Daniel E. Cohen
\endauthor
\affil Queen Mary and Westfield College, London University
\endaffil
\address School of Mathematical Sciences, Queen Mary and Westfield College, 
Mile End Road, London E1 4NS, England
\endaddress
\email D.E.Cohen\@uk.ac.qmw.maths
\endemail

\keywords graph products, isoperimetric functions, Thue systems
\endkeywords
\abstract Let $\Gamma$ be a finite graph, and for each vertex $i$ let $G_i$ 
be a finitely presented group. Let $G$ be the graph product of the $G_i$. That 
is, $G$ is the group obtained from the free product of the $G_i$ by factoring 
out by the smallest normal subgroup containing all $[g,h]$ where $g\in G_i$ 
and $h\in G_j$ and there is an edge joining {\it i\/} and {\it j}. We show 
that $G$ has an isoperimetric function of degree $k>1$ (or an exponential isoperimetric 
function) if each vertex group has such an isoperimetric function.
\endabstract
\endtopmatter
\document \head Graph Products \endhead

Let $\Gamma$ be a finite graph; that is, $\Gamma$ consists of a finite set 
of vertices and a finite set of edges, where each edge is an unordered pair 
of vertices. Let us be given a group $G_i$ for each vertex $i$. Then the {\it 
graph product G\/} of the $G_i$ is the group obtained from the free product 
of the $G_i$ by factoring out by the smallest normal subgroup containing all 
$[g,h]$ where $g\in G_i$ and $h\in G_j$ and there is an edge joining {\it i\/} 
and {\it j}. Note that if $G_i$ has presentation $\langle A_i;R_i\rangle$, 
where the $A_i$ are disjoint, then {\it G\/} has presentation $\langle \bigcup 
A_i;\bigcup R_i \cup S \rangle$, where {\it S\/} is the set of commutators 
$[a,b]$ for $a \in A_i, b \in A_j$, where  {\it i\/} and {\it j\/} are joined 
by an edge. The free product and the direct product are examples of graph products 
(corresponding to graphs with no edges and complete graphs, respectively). 
All groups considered will be finitely presented.

Gersten [G] defines an {\it isoperimetric function} for a finite presenta- 
tion $\langle Y;S \rangle$ of a group $H$ to be a function $f$ such that if 
$w$ is a word of length $n$ in the free group on $Y$ and $w$ equals 1 in $H$ 
then $w$ is the product of at most $f(n)$ conjugates of elements of $S$ and 
their inverses. He shows that if we change to another finite pre- sentation 
then there are positive constants $a, \ldots ,e$ such that the new presentation 
has an isoperimetric function $g$ given by $g(n) = af(bn+c)+dn+e$.

Consequently, we say that $g \preceq f$ if there are positive constants $a, 
\ldots ,e$ such that $g(n) \leq af(bn+c)+dn+e$ for all$n$, and we call $g$ 
equivalent to $f$ is we have both $g \preceq f$ and $f \preceq g$. This is 
slightly different from Gersten's definition of equivalence of functions. I 
prefer this definition because it makes all polynomials of a given degree equivalent, 
and also makes all exponentials equivalent.

When the free monoid $Y^*$ maps onto $H$ (and not just the free group on $Y$) 
we say that $Y$ is a set of {\it monoid generators} of $H$. It is particularly 
useful if $Y$ has the property that to each $y \in Y$ there is $\bar y \in 
Y$ such that $y \bar y$ equals 1 in $H$. When this happens, it is easy to see 
that we can find a set of defining relators containing all the elements $y 
\bar y$ and lying in $Y^*$. It is also easy to check that any finite presentation 
can be changed to a finite presentation of this sort, and that, in looking 
for an isoperimetric function, we need only consider elements of $Y$ and not 
general elements of the free group on $Y$. In this paper we prove the following 
theorem.

\proclaim{Theorem} If each vertex group has an isoperimetric function which 
is polynomial of degree $k>1$ (or an exponential isoperimetric function) then 
so does their graph product.\endproclaim

The theorem will also hold for other classes of isoperimetric functions (this 
follows immediately from the proof), but the precise condition is messy and 
these two cases are the most important. One requirement is that the function 
is at least quadratic. When this holds, it is sufficient that the equivalence 
class contains a function $f$ such that $f(m+n) \leq f(m) + f(n)$ for all $m$ 
and $n$. Ol'shanskii has proved [O] that groups whose isoperimetric function 
is subquadratic are hyperbolic and hence have linear isoperimetric function. 
Note that the graph product of groups with a linear isoperimetric function 
usually does not have a linear isoperimetric function. The inspiration for 
this paper came from work on graph products by Hermiller and Meier [HM]. Her 
discussion of normal forms in graph products, and a similar discussion by Laurence 
[L], led me to the approach given here.

 In proving the theorem we may take any finite presentations of the vertex 
groups. It will be convenient to take the $A_i$ to be disjoint finite sets 
which are monoid generators of $G_i$, so that there is a homomorphism $\pi 
_i:A_i^* \to G_i$ (where, for any set $Y$, $Y^*$ is the free monoid on $Y$).

A non-trivial element of $A_i^*$ will be called an {\it i\/}-word. To each 
{\it i\/}-word we take a symbol $[u]$. Let {\it X\/} be the set of all such 
symbols. Then there is a homomorphism $\rho$ from $X^*$ onto {\it G\/} which 
sends $[u]$ to $\pi _i u$ when {\it u\/} is an {\it i\/}-word. An element of 
$X^*$ will just be called a {\it word}. We say that $[u]$ is {\it in the star 
of i\/} if $u$ is a {\it j\/}-word where $i$ and $j$ are joined by an edge. 
We say that the word $W$ is in the star of $i$ if $W$ is $[u_1]\ldots [u_m]$ 
with each of $[u_1],\ldots ,[u_m]$ in the star of {\it i}.

A sequence of words $W_1,\ldots ,W_n=\epsilon$, where $\epsilon$ is the empty 
word, will be called a {\it reduction sequence\/} if, for all $m<n$, $W_{m+1}$ 
is obtained from $W_m$ by one of the following moves: \roster
\item replace $P[u][v]Q$ by $P[uv]Q$, for any words $P,Q$ and, for any {\it 
i}, any {\it i\/}-words $u$ and $v$;
\item replace $P[u]Q[v]T$ by $P[uv]QT$ for any words $P,T$, any {\it i\/}-words 
$u,v$, any {\it i}, and any word $Q$ in the star of {\it i};
\item replace $P[u]Q$ by $PQ$ for any words $P,Q$ any {\it i}, and any {\it 
i\/}-word {\it u\/} such that $\pi _i u=1$.
\endroster

We refer to {\it i-moves} if there is a need to mention {\it i\/} explicitly.

The following lemma will be proved in the next section.

\proclaim{Lemma}If $\rho W=1$ then there is a reduction sequence starting with 
{\it W}.
\endproclaim

Let $W=W_1,\ldots ,W_n=\epsilon$ be a reduction sequence. We show how to replace 
it by another reduction sequence with nice properties.

Since the sequence ends with $\epsilon$, a move of type 3 must be used at some 
point. Let the first such move be an {\it i\/}-move, going from $W_m$ to $W_{m+1}$. 
Since all earlier moves are of types 1 and 2, it is easy to check that, in 
the sequence $W_1,\ldots ,W_{m+1}$, a {\it j\/}-move followed by an {\it i\/}-move 
can be replaced by an {\it i\/}-move followed by a {\it j\/}-move. Thus we 
may assume that each of the first $m$ moves is an {\it i\/}-move.

We can now see easily (by induction, looking at the reduction sequence beginning 
with $W_2$) that $W$ must be of the form $P[u_1]Q_1[u_2]\ldots Q_{r-1}[u_r]P'$, 
where $u_1,\ldots ,u_r$ are {\it i\/}-words, $Q_1,\ldots ,Q_{r-1}$ are (possibly 
empty) words in the star of {\it i}, and $\pi _i(u_1\ldots u_r)=1$.

For an arbitrary word $V=[v_1]\ldots [v_s]$, define $\|V\|$ to be $|v_1|+\cdots 
+|v_s|$, where $|v|$ is the length of $v$. We define the {\it weight\/} of 
a move of type 1 to be 0, the weight of a move of type 2 to be $\|Q\|\cdot 
|v|$, and the weight of a move of type 3 to be $f(|u|)$, where $f$ is an isoperimetric 
function for all of the groups $G_i$. We define the weight of a reduction sequence 
to be the sum of the weights of its moves, and we define the weight of a word 
$W$ for which $\rho W=1$ to be the minimum weight of the reduction sequences 
beginning with $W$.

Let $g(n)$ be the maximum of $f(n_1)+\cdots +f(n_s)$ over all $s$ and all positive 
integers $n_1, \ldots ,n_s$ whose sum is $n$. Note that if $f$ is polynomial 
of degree $k$ (or exponential) then so is $g$.

We next show that the weight of a word $W$ with $\rho W=1$ is at most $\|W\|^2+g(\|W\|)$. 
As already remarked, we can write $W$ as $P[u_1]Q_1[u_2]\ldots Q_{r-1}[u_r]P'$, 
where $u_1,\ldots ,u_r$ are {\it i\/}-words, $Q_1,\ldots ,Q_{r-1}$ are (possibly 
empty) words in the star of {\it i}, and $\pi _i(u_1\ldots u_r)=1$. Then there 
is a reduction sequence beginning with $W,P[u_1u_2]Q_1Q_2[u_3]\ldots [u_r]P',\ldots 
P[u_1\ldots u_r]Q_1 \ldots Q_{r-1},W'$, where $W'$ is $PQ_1 \ldots Q_{r-1}P'$. 
Since the sum of the weights of the moves from $W$ to $W'$ is at most $ \|W\| 
\cdot (|u_1|+ \cdots +|u_r| + f(|u_1 \ldots u_r|)$, the required result holds 
by induction.

Finally, we use this result to obtain an isoperimetric function for $G$. We 
use the set of monoid generators $\bigcup A_i$. There are homomorphisms $\pi 
:(\bigcup A_i)^* \to G$, $\alpha :(\bigcup A_i)^* \to X^*$, and $\beta :X* 
\to (\bigcup A_i)^*$, defined by $\pi a = \pi _ia$ for $a \in A_i$, $\alpha 
a = [a]$, and $\beta [u] =u$. Plainly, for any $w \in (\bigcup A_i)^*$, we 
have $\beta \alpha w =w$ and $\rho \alpha w = \pi w$. Also $|w|=\|\alpha w\|$.

It is easy to see that if $W'$ is obtained from $W$ by a move of weight $k$ 
then $\beta W$ is the product of $\beta W'$ and $k$ conjugates of the defining 
relators (and their inverses) of the finite presentation of $G$. By induction 
on the length of the reduction sequence, if $\rho W = 1$ then $\beta W$ is 
the product of at most weight\,($W$) conjugates of the defining relators and 
their inverses.

Applying this to $\alpha w$, where $\pi w=1$, and using the formula for the 
weight, we see that $g(n)+n^2$ is an isoperimetric function for our presentation 
of $G$, proving the theorem.

\head Thue Systems \endhead

Let $X$ be an arbitrary set. A {\it Thue system\/}, or {\it rewriting system\/} 
on $X$ is a subset $S$ of $X^* \times X^*$. Such a system induces an equivalence 
relation on $X^*$; namely, the smallest equivalence relation such that $ulv 
\equiv urv$ for all words $u$,$v$ and all pairs $(l,r)$ in $S$. The quotient 
of $X^*$ by this equivalence is called {\it the monoid presented by $\langle 
X;S \rangle $}.

When we look for normal forms for the equivalence classes, there are two ways 
to proceed. One treats all members of $S$ alike, and compares the equivalence 
relation with the  non-symmetric relation in which we can replace $ulv$ by 
$urv$ but not vice versa. We then endeavour to see if this terminates, and 
whether it provides a unique normal form.

The other approach, which is more convenient in our situation, begins by assuming 
that $S$ consists only of pairs for which $|l| \ge |r|$ and that, for any $(l,r)$ 
in $S$ with $|l|=|r|$ we also have $(r,l)$ in $S$. This can be done without 
loss of generality, since we get the same equivalence relation if we replace 
a pair $(l,r)$ by $(r,l)$ and also if we add pairs $(r,l)$ for which $(l,r)$ 
are already in $S$.

If we do this, then, when we consider replacing $ulv$ by $urv$ but not vice 
versa, if $|l|=|r|$ we can use the further pair $(r,l)$ to return from $urv$ 
to $ulv$. Consequently, we treat such pairs differently from thos pairs for 
which $|l|>|r|$. It is quite common in computer science to distinguish between 
the two approaches by using the phrase `rewriting system' for the first one 
and the phrase `Thue system' for the second.

We write $ulv \to urv$ for a pair $(l,r)$ with $|l|>|r|$, and $ulv \sim urv$ 
for a pair with $|l| = |r|$. We let $\overset *\to\rightarrow$ and $\overset 
*\to\sim$ be the reflexive transitive closures of these.

We say that the pair $u,v$ is {|it almost confluent\/} if there are 
$u_1,v_1$ such that $u \overset *\to\rightarrow u_1$,$v \overset 
*\to\rightarrow v_1$, and  $u_1 \overset *\to\sim v_1$. Plainly, almost 
confluent words are equivalent, and we call $S$ almost confluent if every 
pair of equivalent words is almost confluent.
In searching for nice representatives of the equivalence classes, it is 
particularly helpful if $S$ is almost confluent.
Clearly, when this holds, if $u$ is equivalent to $\epsilon $ then $u 
\overset *\to\rightarrow \epsilon$.

This situation is very familiar to computer scientists. A sufficient 
condition for the property to hold can be given in terms of the behaviour 
of certain {\it critcal pairs\/} of words, which arise from certain words 
in which two of the elements of $S$ may be used.  The situation is less 
well-known to group theorists, but results sometimes referred to as `Peak 
Reduction Lemmas' are essentially of this form.

Huet [H] showed that $S$ is almost confluent whenever almost confluence 
holds for all pairs $u,v$ such that there is some $w$ with $w \rightarrow 
u$ and either $w \rightarrow v$ or $w \sim v$. It is not difficult to 
prove this directly using peak reduction arguments. (If readers want to 
look at [H], they should note that Huet's $\sim$ is our $\overset
*\to\sim$.)

Huet also showed that we do not even have to consider all such pairs. It 
is enough to consider those $w$ of form $abc$ for some words $a,b,c$ such 
that $S$ either has elements $(ab,r_1)$ and $(bc,r_2)$ or has elements 
$(abc,r_1)$ and $(b,R_2)$, with $u$ and $v$ (or $v$ and $u$) obtained from 
$w$ by applying these two elements. These pairs $u,v$ of words are called 
the critical pairs. Huet's proof applies in much more generality, and it 
is probably simpler to prove this directly in our situation.

We now return to graph products, with the set $X$ as in the first section. 
We shall prove the lemma by applying this theory of Thue systems. The set $S$ 
will consist of the following pairs:

\roster
\item $([u][v],[uv])$, where $u$ and $v$ are {\it i\/}-words for some $i$;
\item $([u]P[v],[uv]P)$, where $u$ and $v$ are {\it i\/}-words for some $i$ 
and $P$ is in the star of $i$;
\item $([u],\epsilon )$, where $u$ is an {\it i\/}-word such that $\pi _i u=1$;
\item $([u],[v])$,  where $u$ and $v$ are {\it i\/}-words for some $i$ and 
$\pi _i u = \pi_i v$;
\item $([u][v],[v][u])$, where $u$ is an {\it i\/}-word, $v$ is a {\it j\/}-word, 
and there is an edge of the graph joining $i$ and $j$.
\endroster

It is clear that the set of {\it i\/}-words for a given $i$, together with 
the corresponding pairs of types 1,3, and 4, form a monoid presentation for 
$G_i$; this is just a variant of the multiplication table presentation. If 
we take all pairs of types 1,3,4, and 5 we then clearly obtain a monoid presentation 
for $G$. We can then add the pairs of type 2 and still get a monoid presentation 
for $G$, since the two elements of a pair of type 2 clearly give the same element 
of $G$.

To prove the lemma, we need only show that the criterion mentioned above 
is satisfied.

First look at $[u]P[v]Q[w]$, where $u$, $v$, and $w$ are {\it i\/}-words for 
some $i$, and $P$ and $Q$ are (possibly empty) words in the star of $i$. We 
have $[u]P[v]Q[w] \to [uv]PQ[w]$ and also $[u]P[v]Q[w] \to [u]P[vw]Q$. Here 
we find that $[uv]PQ[w] \to [uvw]PQ$ and also $[u]P[vw]Q \to [uvw]PQ$.

Next, look at $[u]P[v]$, where $u$ and $v$ are {\it i\/}-words, $P$ is in the 
star of $i$, and $\pi _i v = 1$. Then $[u]P[v] \to [uv]P$ and also $[u]P[v] 
\to [u]P$. Since $\pi _i (uv) = \pi _i u$, we have $[uv]P \sim [u]P$, using 
a pair of type 4. If we have $\pi _i u =1$ instead of $\pi _i v =1$, then $[u]P[v] 
\to [uv]P$ and $[u]P[v] \to P[v]$. Here we have $[uv]P \sim [v]P$, using a 
pair of type 4, and $[v]P \overset\*\to\sim P[v]$, using pairs of type 5 (since 
$P$ is in the star of $i$).

Suppose we have a word $[u]P[v]$, where $u$ and $v$ are {\it i\/}-words for 
some $i$ and $P$ is a (possibly empty) word in the star of $i$. Let $w$ be 
an {\it i\/}-word such that $\pi _i w = \pi _i u$. Then $[u]P[v] \to [uv]P$ 
and $[u]P[v] \sim [w]P[v]$. We then have $[w]P[v] \to [wv]P$ and $[uv]P \sim 
[wv]P$. A similar argument works when, instead of $\pi _i w = \pi _i u$, we 
have $\pi _i w = \pi _i v$.

Let $u$ be an {\it i\/}-word such that $\pi _i u =1$, and let $w$ be an {\it 
i\/}-word such that $\pi _i w = \pi _i u$. Then $[u] \to \epsilon$ and $[u] 
\sim [w]$, and we also have $[w] \to \epsilon$. Let $v$ be a {\it j\/}-word, 
where is an edge joining $i$ and $j$. Then we have $[u][v] \sim [v][u]$ and 
also $[u][v] \to [v]$. Here we have $[v][u] \to [v]$.

Suppose we have a word $[u]P[v]$, where $u$ and $v$ are {\it i\/}-words for 
some $i$, and $P$ is a (possibly empty) word in the star of $i$. Let $w$ be 
a {\it j\/}-word, where there is an edge joining $i$ and $j$. Then $[u]P[v][w] 
\to [uv]P[w]$ and $[u]P[v][w] \sim [u]P[w][v]$. Since $P[w]$ is in the star 
of $i$, we have $[u]P[w][v] \to [uv]P[w]$. Finally, we have $[w][u]P[v] \to 
[w][uv]P$ and $[w][u]P[v] \sim [u][w]P[v]$. We then have $[u][w]P[v] \to [uv][w]P$ 
and $[w][uv]P \sim [uv][w]P$.

We have now shown that all critical pairs satisfy the required criterion, and 
the lemma is proved.


\ref \key G \by S.M. Gersten \paper Isoperimetric and Isodiametric 
Functions of Finite Presentations \inbook Proceedings of the Conference on 
Geometric Group Theory (Isle of Thorns, 1992) \endref
\ref \key HM \by S. Hermiller and J. Meier \paper Algorithms and Geometry for 
Graph Products of Groups \toappear \endref
\ref \key H \by G. Huet \paper Confluent Reduction: Abstract Properties 
and Applications to Term Rewriting Systems \jour J. Ass. Computing 
Machinery \vol 27 \yr1980 \pages 797--821 \endref
\ref \key L \by M. R. Laurence \book Automorphisms of graph products of groups 
\bookinfo Ph.D. Thesis \publ Queen Mary and Westfield College \publaddr London 
\yr1993 \endref

\ref \key O \by A. Y. Ol'shanskii \paper Hyperbolicity of groups with subquadratic 
isoperimetric inequality \yr1991 \vol 1 \pages 281 -- 289 \jour Int. J. Alg. 
Comp. \endref
\enddocument